----------
X-Sun-Data-Type: default
X-Sun-Data-Description: default
X-Sun-Data-Name: uncertainty
X-Sun-Charset: us-ascii
X-Sun-Content-Lines: 624

%Beginning Formats                                        

\magnification=1200
%\magnification=1000
%\magnification=\magstep1
%\vsize=20truecm
%\voffset=1.00truein
\settabs 18 \columns
%\hoffset=3.75truecm
%\hoffset=1.00truein
%\hsize=14truecm

%\nopagenumbers
\baselineskip=17 pt
%\baselineskip=12pt
%\ifnum\pageno=1
%\topinsert \vskip 1.00 in
%\endinsert
%\vsize=7.5in
%\fi

\def\s{\smallskip}

\def\sqr#1#2{{\vcenter{\vbox{\hrule height.#2pt
 \hbox{\vrule width.#2pt height#1pt \kern#1pt
 \vrule width.#2pt} \hrule height.#2pt}}}}

\def\operp{\hbox{${\kern+.25em{\bigcirc}
\kern-.85em\bot\kern+.85em\kern-.25em}$}}
%\def\gapprox
%{\hbox{$
%\smash{lower0.5ex\hbox{$\scriptstyle>$}} \atop
%\smash{raise0.3ex\hbox{$\scriptstyle \sim$}}
%$}}
%\def\lapprox
%{\hbox{$
%\smash{lower0.5ex\hbox{$\scriptstyle<$}} \atop

%\smash{raise0.3ex\hbox{$\scriptstyle \sim$}}
%$}}
\def\lsim{\;\raise0.3ex\hbox{$<$\kern-0.75em\raise-1.1ex\hbox{$\sim$}}\;}
\def\gsim{\;\raise0.3ex\hbox{$>$\kern-0.75em\raise-1.1ex\hbox{$\sim$}}\;}
\def\no{\noindent}

\def\ce{\centerline}
\def\ve{\vfill\eject}
\def\rdots{\mathinner{\mkern1mu\raise1pt\vbox{\kern7pt\hbox{.}}\mkern2mu
 \raise4pt\hbox{.}\mkern2mu\raise7pt\hbox{.}\mkern1mu}}

\def\e e{$e^+ e^-$ }

%End of Beginning Formats
%Beginning of Letter Heading

\rightline{UCLA/97/TEP/14}
\rightline{May 1997}
\vskip2.0cm

\ce{\bf $q$--UNCERTAINTY RELATIONS}
\vskip.5cm
\ce{Robert J. Finkelstein}
\s
\ce{\it Department of Physics and Astronomy}
\ce{\it University of California, Los Angeles, CA 90095-1547}
\vskip1.0cm

\no {\bf Abstract.}  As one would anticipate from the realization
of the $q$-commutators by difference operators, the states of maximum
localization are smeared $\delta$-functions depending on $q$.
\vskip.5cm

\line{{\bf 1. Introduction.} \hfil}
\s

The Heisenberg uncertainty rules, depending on the standard commutators
between conjugate observables, will of course be modified if the
commutators are deformed.  This question, which has been previously
studied by A. Kempf,$^1$ will be further examined here and illustrated
for the $q$-harmonic oscillator.
\vskip.5cm

\line{{\bf 2. $q$-Oscillator.} \hfil}
\s

In the Fock representation the familiar equations of the $q$-oscillator are
$$
\eqalignno{& a\bar a-q\bar aa = \Delta & (2.1) \cr
& H = {1\over 2}(a\bar a+\bar aa) & (2.2) \cr}
$$
\no where we consider two realizations of (2.1) which we shall express
in terms of eigenstates $|n\rangle$ of $H$:
$$
\eqalignno{{\hbox{(a)}}\quad \quad \quad \quad 
& a|n\rangle = \langle n\rangle^{1/2}
|n-1\rangle & (2.3) \cr
~~~ \quad \quad \quad \quad & \bar a|n\rangle 
= \langle n+1\rangle^{1/2}|n+1\rangle & (2.4) \cr
~~~ \quad \quad \quad \quad & H|n\rangle = {1\over 2}\bigl(\langle n\rangle 
+ \langle n+1\rangle\bigr)
& (2.5) \cr
~~~ \quad \quad \quad \quad & \Delta |n\rangle = |n\rangle & (2.6) \cr}
$$
$$
\eqalignno{{\hbox{(b)}} \quad \quad \quad \quad
& a|n\rangle = [n]^{1/2}|n-1\rangle & (2.7) \cr
~~~ \quad \quad \quad \quad &\bar a|n\rangle 
= [n+1]^{1/2}|n+1\rangle & (2.8) \cr
~~~ \quad \quad \quad \quad & H|n\rangle = {1\over 2}([n]+[n+1]) 
& (2.9) \cr
~~~ \quad \quad \quad \quad & \Delta |n\rangle = q^{-n}|n\rangle & (2.10) \cr}
$$
\no where
$$
\eqalignno{\langle n\rangle &= {q^n-1\over q-1} & (2.11) \cr
[n] &= {q^n-q^{-n}\over q-q^{-1}} & (2.12) \cr}
$$

Set $A = \left(\matrix{a \cr \bar a \cr}\right)~.$  Then Eq. (2.1) may be
written as either
$$
A^t\epsilon A = q^{-1/2} \quad~~ \hbox{if} \quad \Delta =1  \eqno(2.13)
$$
\no or
$$
A^t\epsilon A = q^{-\bigl({n+1\over 2}\bigr)} \quad \hbox{if} \quad
\Delta = q^{-n}  \eqno(2.14)
$$
\no where
$$
\epsilon = \left(\matrix{0 & q^{-1/2} \cr
-q^{1/2} & 0 \cr}\right)~. \eqno(2.15)
$$

\no The $\epsilon$ matrix is invariant under transformation $T$ belonging
to $SU_q(2)$:
$$
T^t\epsilon T = T\epsilon T^t = \epsilon~. \eqno(2.16)
$$
\no Therefore any new vector $X$ defined by
$$
A = TX \eqno(2.17)
$$
\no will satisfy
$$
X^t\epsilon X = \Delta~. \eqno(2.18)
$$

The transnformation from Fock observables to configuration space may be
accomplished by choosing
$$
X = \left(\matrix{D \cr x\cr}\right) \quad p = {\hbar\over i} ~D
\eqno(2.19)
$$
\no Then
$$
qxp-px = i\hbar\Delta~. \eqno(2.20)
$$
\no If $\Delta = 1$, then (2.20) requires
$$
D^q = {1\over x} \langle \theta\rangle_q = {1\over x}
{q^\theta-1\over q-1}~. \eqno(2.21)
$$
\no Here $\theta$ is the dilatation operator
$$
\theta = x{d\over dx}~. \eqno(2.22)
$$
\no If $\Delta = q^{-n}$
$$
D = {1\over x} [\theta]_q = {1\over x}
{q^\theta-q^{-\theta}\over q-q^{-1}}~. \eqno(2.23)
$$
\no By (2.21)
$$
D^q f(x) = {f(qx)-f(x)\over (q-1)x} ~. \eqno(2.24)
$$
\no By (2.23)
$$
Df(x) = {f(qx)-f(q_1x)\over (q-q_1)~x} \eqno(2.25)
$$
\no where
$$
q_1=q^{-1}~.
$$

Although the difference operators $D^q$ and $D$ convey the same idea, it
turns out that $D$ is the proper choice in constructing the Fourier
transform between configuration and momentum space.  (In previous work$^2$
we used the form (2.24) corresponding to $\langle~~\rangle$).  We shall
now use (2.25) or $[~~]$.  Then
$$
\eqalignno{a &= \alpha D + \beta x & (2.26) \cr
\bar a &= -q_1\beta D + \bar\alpha x & (2.27) \cr}
$$
\no where
$$
T = \left(\matrix{\alpha & \beta \cr -q_1\bar\beta & \bar\alpha \cr}\right)
\eqno(2.28)
$$
\no and $D$ is given by (2.23).
\vskip.5cm

\line{{\bf 3. Configuration Space Representation of Fock States.} \hfil}
\s

These are obtained by the Fock ladder operators.$^2$  The ground state was
previously determined by
$$
a^q|0\rangle = (\beta x+\alpha D^q)|0\rangle = 0 \eqno(3.1)
$$
\no and the higher states by the raising operators.  We illustrate the
different procedures required by $D^q$ and $D$ by displaying the ground
states in both cases.  The solution of (3.1) for 
$\langle x|0\rangle = \psi_o^q(x)$ is
$$
\psi_o^q(x) = \prod^\infty_o (1+Kq^{2s})^{-1} \psi^q(0)~,
\quad q< 1 \eqno(3.2a)
$$
\no where
$$
K = (1-q)\alpha^{-1}\beta x^2~. \eqno(3.2b)
$$
\no Let us now obtain the ground state $\psi_o(x)$ corresponding to $D$
by solving
$$
(\beta x+\alpha D)\psi_o(x) = 0 \eqno(3.3)
$$
\no where $D$ is given by (2.23).  Then
$$
\psi_o(qx) - \psi_o(q_1x) = K^\prime x^2\psi_o(x) \eqno(3.4)
$$
\no where
$$
K^\prime = (q_1-q) \alpha^{-1}\beta~. \eqno(3.5)
$$
\no The difference equation (3.4) relates $\psi_o(x)$ at three distinct
points.  To reduce this equation to the previous case where there are
only two distinct points define
$$
\varphi(qx) = \psi_o(qx)/\psi_o(x)~. \eqno(3.6)
$$
\no Then by (3.4)
$$
\varphi(qx)-\varphi(q_1x) = K^\prime x^2 \eqno(3.7)
$$
\no or
$$
\varphi(x) = \varphi(q^2x)-K^\prime q^2x^2~. \eqno(3.8)
$$
\no By iterating this relation one obtains
$$
\varphi(qx) = \varphi(q^{2n+1}x)-K^\prime q^4x^2
{q^{4n}-1\over q^4-1}~. \eqno(3.9)
$$
\no Letting $n\to\infty$ one finds
$$
\varphi(qx) = \varphi(0) - {K^\prime q^4x^2\over 1-q^4} \eqno(3.10)
$$
\no since $q<1$.  As $\varphi(0)=1$, one may by (3.6) rewrite (3.10)
in the following form
$$
\psi_o(qx) = \biggl[1-{K^\prime q^4x^2\over 1-q^4}\biggr]\psi_o(x)~. \eqno(3.11)
$$
\no Again by iteration $\psi_o(x)$ may be found as the infinite product
$$
\psi_o(x) = \prod^\infty_o \bigl(1-Lq^{2s}\bigr)^{-1}
\psi_o(0) \eqno(3.12)
$$
\no where
$$
\eqalign{L &= {K^\prime q^4x^2\over 1-q^4} \cr
&= {q^3x^2\over 1+q^2} \alpha^{-1}\beta~. \cr} \eqno(3.13)
$$
\no Both (3.2) and (3.12) are formally convergent for all $x$.  Since these
series actually lie in the $SU_q(2)$ algebra, however, their interpretation
as state functions must depend on supplementary rules for their numerical
valuation, as discussed elsewhere.$^{2,3}$

By applying the raising operator the complete set of excited states may be
found.  In a similar way the complete set of momentum states may be found.
The ground states are $q$-Gaussian and all states vanish at $x=\infty$.$^2$
\vskip.5cm

\line{{\bf 4. Alternative Transition to Configuration Space.} \hfil}
\s

If the position and momentum variables are defined by the standard
expressions$^1$ one has
$$
\eqalign{x &= L(\bar a+a) \cr
p &= iK(\bar a-a)~.} \eqno(4.1)
$$
\no If $a$ and $\bar a$ were hermitian conjugate then $x$ and $p$ would be
separately hermitian.  Since $q \not= 1$, however, $a$ and $\bar a$ are not
hermitian conjugate, and $x$ and $p$ are not hermitian.  The transformation
(4.1) does not preserve the $q$-commutator.  Theh usual
commutator is$^1$
$$
(x,p) = i\hbar\biggl[\Delta + (q-1){1\over 4}
\biggl({x^2\over L^2} + {p^2\over K^2}\biggr)\biggr] \eqno(4.2)
$$
\no if one sets
$$
LK = {\hbar\over 2}{q+1\over 2}~. \eqno(4.3)
$$
\no The corresponding uncertainty relations are$^1$
$$
\bigl(\langle(\Delta x)^2\rangle\langle(\Delta p)^2\rangle\bigr)^{1/2}\geq
{\hbar^2\over 4}\biggl[\Delta + (q-1)
\biggl({\langle x^2\rangle\over 4L^2}+{\langle p^2\rangle\over 4K^2}
\biggr)\biggr]~. \eqno(4.4)
$$
\no Since $x$ and $p$ are not the usual hermitian observables supplementary
rules of interpretation would also be necessary here.

The Fock commutators (2.1) imply the failure of translational invariance
and the modified uncertainty relations as exemplified in (4.4), if the
transition to configuration space is made by the standard relations (4.1).
Failure of translational invariance and deformed uncertainty relations
must always appear, however, in whatever way the transition from (2.1) to
configuration space is made.  Let us therefore examine this question if
$x$ and $p$ are defined by (2.26) and (2.27).  To do this we shall next
introduce a $q$-Fourier transform.
\vskip.5cm

\line{{\bf 5. $q$-Fourier Transforms.} \hfil}
\s

One may define ``eigenstates" of momentum by
$$
{\hbar\over i} D_x\psi_p(x) = p\psi_p(x)~. \eqno(5.1)
$$
\no Then
$$
\psi_p(x) = \sum {(ipx/\hbar)^n\over [n]!}~. \eqno(5.2)
$$
\no If $D_x^q$ instead of $D_x$ had apppeared in (5.1), then $\langle n\rangle$
would have replaced $[n]$ in (5.2).  The series (5.2) is convergent for
large $x$.  The corresponding $\langle n\rangle$ series is not convergent
for large $x$; for this reason the symmetric derivative is required.

In standard quantum mechanics one transforms from configuration to momentum
space by writing the wave function in $x$-space as a sum over momentum
eigentates.  The corresponding procedure here is
$$
\psi(x) = \sum_p \psi_p(x) \varphi(p)~. \eqno(5.3)
$$
\no In Dirac notation
$$
\psi_p(x) = \langle x|p\rangle \eqno(5.4)
$$
\no and (5.3) becomes
$$
\langle x|~\rangle = \sum_p\langle x|p\rangle\langle p|~\rangle~. \eqno(5.5)
$$
Here we choose $\Sigma_p$ to be the $q$-integral.  Then
$$
\psi(x) = \int^\infty_{-\infty} {\cal{E}}(ipx)\varphi(p) d_qp \eqno(5.6)
$$
\no where we have set
$$
\langle x|p\rangle = {\cal{E}}(ipx)~. \eqno(5.7)
$$

Then (5.6) is the $q$-Fourier decomposition of $\psi(x)$.  In order that
the series for the kernel ${\cal{E}}(ipx)$ be convergent for all $x$ it
is necessary to use (5.2) instead of the $\langle n\rangle$ series.

To do integration by parts we need the Leibniz rule for $D$, namely:
$$
\eqalign{D_x[f(x)g(x)] &= f(qx)D_xg(x) + g(q_1x)D_xf(x) \cr
\noalign{\hbox{or}}
&= f(q_1x) D_xg(x) + g(qx) D_xf(x)~. \cr} \eqno(5.8)
$$
\no Let us
 test (5.6) by computing
$$
\eqalignno{a_x\psi &= (\alpha D_x+\beta x)\int^\infty_{-\infty}
{\cal{E}}(ipx) \varphi(p) d_qp & (5.9) \cr
&= \int^\infty_{-\infty} (\alpha(ip)+\beta x)
{\cal{E}}(ipx)\varphi(p) d_qp & (5.10) \cr}
$$
\no where we have used
$$
D_x{\cal{E}}(ipx) = ip{\cal{E}}(ipx)~. \eqno(5.11)
$$
\no Note
$$
\eqalignno{\int^\infty_{-\infty} D_p\bigl({\cal{E}}(ixp)\varphi(p)\bigr)
d_qp &= \bigl({\cal{E}}(ixp)\varphi(p)\bigr)^\infty_{-\infty} & (5.12)\cr
&= 0 & (5.13) \cr}
$$
\no since $\int_q$ is inverse to $D$ and since
$$
\varphi(\infty) = \varphi(-\infty) = 0~. \eqno(5.14)
$$
\no By (5.13) and (5.8)
$$
\int^\infty_{-\infty}\bigl({\cal{E}}(iq_1px)D_p\varphi(p) +
\varphi(qp) D_p{\cal{E}}(ipx)\bigr) d_qp = 0~. \eqno(5.15)
$$

Corresponding to (5.11) we have
$$
D_p{\cal{E}}(ipx) = ix{\cal{E}}(ipx)~. \eqno(5.16)
$$
\no By (5.15) and (5.16) we have
$$
\int^\infty_{-\infty} \varphi(qp)ix{\cal{E}}(ipx)d_qp =
-\int^\infty_{-\infty} {\cal{E}}(iq_1px)D_p\varphi(p) d_qp~. \eqno(5.17)
$$
\no Set
$$
\eqalign{p^\prime &= qp \cr
x^\prime &= q_1x~.\cr} \eqno(5.18) 
$$
\no Then
$$
\eqalignno{\int^\infty_{-\infty} \varphi(p^\prime)
ix^\prime{\cal{E}}(ip^\prime x^\prime) d_qp^\prime &=
-\int^\infty_{-\infty} {\cal{E}}(ipx^\prime) D_p\varphi(p) d_qp & (5.19) \cr
\noalign{\hbox{or}}
\int^\infty_{-\infty}\varphi(p) ix{\cal{E}}(ipx)d_qp &=
-\int^\infty_{-\infty}{\cal{E}}(ipx)D_p\varphi(p) d_qp~. & (5.20) \cr}
$$

By (5.10) and (5.20)
$$
\eqalignno{a_x\psi &= \int^\infty_{-\infty} {\cal{E}}(ipx)
i(\alpha p+\beta D_p)\varphi d_qp & (5.21) \cr
\noalign{\hbox{or}}
a_x\psi &= \int^\infty_{-\infty} {\cal{E}}(ipx) a_p\varphi d_qp & (5.22) \cr}
$$
\no where
$$
a_p = i(\alpha p + \beta D_p)~. \eqno(5.23)
$$

Let us next check the following ansatz for the reverse transformation
$$
\varphi(p) = \int^\infty_{-\infty}{\cal{E}}(-ipx) \psi(x) d_qx~. \eqno(5.24)
$$
\no Then
$$
a_p\varphi(p) = \int^\infty_{-\infty}(i\alpha p + \beta x)
{\cal{E}}(-ippx)\psi(x) d_pK~. \eqno(5.25)
$$
\no Again the right-hand side may be transformed to give
$$
a_p\varphi(p) = \int^\infty_{-\infty}{\cal{E}}(-ipx) a_x\psi(x) d_qx~.
\eqno(5.26)
$$
\no By (5.22) and (5.26) the ground states in $x$ and $p$ space correspond.
Thus
$$
a_x\psi_o(x) = 0~~\longleftrightarrow~~a_p\varphi_o(p) = 0~. \eqno(5.27)
$$
\no The argument may be completed with the aid of the raising operators.
\vskip.5cm

\line{{\bf 6. The $q$-Delta Function.} \hfil}
\s

In Dirac notation, Eqs. (5.6) and (5.24) may be written
$$
\eqalignno{\langle x|~\rangle &= \int \langle x|p\rangle d_qp
\langle p|~\rangle & (6.1) \cr
\langle p|~\rangle &= \int \langle p|x\rangle d_qx
\langle x|~\rangle & (6.2) \cr}
$$
\no Combining these two equations we have
$$
\eqalignno{\langle x|~\rangle &= \int \langle x|p\rangle d_qp
\int \langle p|x^\prime\rangle d_qx^\prime\langle x^\prime|~\rangle &
(6.3) \cr
\noalign{\hbox{or}}
\langle x|~\rangle &= \int \delta_q(x,x^\prime) d_qx^\prime
\langle x^\prime|~\rangle & (6.4) \cr}
$$
\no where
$$
\eqalign{\delta_q(x,x^\prime) &= \langle x|x^\prime\rangle_q \cr
&= \int \langle x|p\rangle d_qp\langle p|x^\prime\rangle \cr
&= \int {\cal{E}}(ipx){\cal{E}}(-ipx^\prime) d_qp~. \cr} \eqno(6.5)
$$
\no To evaluate the integral (6.5) note
$$
\int D_p\bigl({\cal{E}}(ipx){\cal{E}}(-ipx^\prime)\bigr) d_qp =
\bigl({\cal{E}}(ipx){\cal{E}}(-ipx^\prime)\bigr)^\infty_{-\infty}~.
\eqno(6.6)
$$
\no The left side of (6.6) is
$$
\eqalign{&\int\bigl({\cal{E}}(ipx)D_p{\cal{E}}(-ipx^\prime) +
{\cal{E}}(-iq_1px^\prime)D_p{\cal{E}}(ipx)\bigr)d_qp \cr
&= \int\bigl({\cal{E}}(iqpx)(-ix^\prime){\cal{E}}(-ipx^\prime) +
{\cal{E}}(-iq_1px^\prime) ix{\cal{E}}(ipx)\bigr)d_qp~. \cr} \eqno(6.7)
$$

Set $p^\prime=qp$ in the first integral on the right side of (6.7).  Then
(6.7) becomes
$$
\eqalign{&\int{\cal{E}}(ip^\prime x)(-ix^\prime)
{\cal{E}}(-ip^\prime q_1x^\prime)
q_1dp^\prime+\int{\cal{E}}(-iq_1px^\prime)ix{\cal{E}}(ipx)d_qp \cr
&= \int {\cal{E}}(ipx){\cal{E}}(-iq_1px^\prime)(-iq_1x^\prime+ix)d_qp~. \cr}
\eqno(6.8)
$$

Set $q_1x^\prime = x^{\prime\prime}$.  Then (6.8) becomes
$$
i(x-x^{\prime\prime}) \int {\cal{E}}(ipx){\cal{E}}(-ipx^{\prime\prime})
d_qp~. \eqno(6.9)
$$
\no Therefore
$$
\int {\cal{E}}(ipx){\cal{E}}(-ipx^{\prime\prime}) d_qp =
{\bigl({\cal{E}}(ipx){\cal{E}}(-iqpx^{\prime\prime})\bigr)^\infty_{-\infty}
\over i(x-x^{\prime\prime})} \eqno(6.10)
$$
\no or
$$
\delta_q(x,x^\prime) = \lim_{P\to\infty}
\biggl({{\cal{E}}(iPx){\cal{E}}(-iqPx^\prime)-{\cal{E}}(-iPx){\cal{E}}
(iqPx^\prime)\over i(x-x^\prime)}\biggr)~. \eqno(6.11)
$$
\no If $q=1$, $\delta_q(x,x^\prime)$ becomes the Dirac $\delta$-function
$$
\delta(x-x^\prime) = 2 \lim_{P\to\infty}
\biggl({\sin P(x-x^\prime)\over x-x^\prime}\biggr)~. \eqno(6.12)
$$
\no In this case $\langle x|x^\prime\rangle$ vanishes unless $x=x^\prime$.
It is also translationally invariant.

If $q\not= 1$, $\delta_q(x,x^\prime)$ is not translationally invariant.
Set $x^\prime=0$.  Then
$$
\langle x|0\rangle_q = \delta_q(x,0) = \lim_{P\to\infty}
P\biggl({\sin_q Px\over xP}\biggr)~. \eqno(6.13)
$$
\no Set $x=0$.  Then
$$
\langle 0|x^\prime\rangle_q = \delta_q(0,x^\prime) =
\lim_{P\to\infty} 
P\biggl({\sin_q(qPx^\prime)\over x^\prime P}\biggr)~. \eqno(6.14)
$$

The state of maximum localization is therefore a smeared $\delta$-function,
the $q$ $\delta$-function.  This result is consistent with the lattice
properties implied by the implementation of the momentum operator by a
difference operator.  The smeared $\delta$-function is also in agreement
with the higher uncertainty bound expressed by Eq. (4.4).  These results
are purely formal, however, and the pair $(x,p)$ defined by either
((2.26), (2.27)) or (4.1) would require additional rules for their physical
interpretation.

I would like to thank Dr. A. Kempf for a stimulating discussion.
\vskip.5cm

\line{{\bf References.} \hfil}
\s

\item{1.} A. Kempf, J. Math. Phys. {\bf 35}. 4483 (1994).
\item{2.} R. Finkelstein and E. Marcus, J. Math. Phys. {\bf 36}, 2652
(1995).
\item{3.} R. Finkelsteinn, J. Math. Phys. {\bf 37}, 2628 (1996);
\item{} A. C. Cadavid and R. Finkelstein, J. Math. Phys. {\bf 37}, 3675
(1996).

\ve

\bye